\documentclass[11pt,letterpaper]{article}

\title{Automated Counting of Restricted Motzkin Paths}
\author{AJ Bu}
\date{\today}

\usepackage[T1]{fontenc}
\usepackage{lmodern}
\usepackage[pdftex,margin=1in]{geometry}
\usepackage{enumerate}
\usepackage{amsmath}
\usepackage{amssymb}
\usepackage{fancyhdr}
\usepackage{lastpage}
\usepackage[alwaysadjust]{paralist}
\usepackage{amsthm}
\usepackage{graphicx}

\usepackage{wasysym}
\usepackage{xcolor}
\usepackage{listings}
\usepackage{hyperref}

\newcommand{\la}{\langle}
\newcommand{\ra}{\rangle}

\newtheorem{Thm}{Theorem}

\theoremstyle{definition}
\newtheorem{Def}[Thm]{Definition}

\theoremstyle{remark}

\theoremstyle{definition}

\theoremstyle{definition}

.

\lstdefinelanguage{Maple}% 
{morekeywords={and,assuming,break,by,catch,description,do,done,% 
		elif,else,end,error,export,fi,finally,for,from,global,if,% 
		implies,in,intersect,local,minus,mod,module,next,not,od,% 
		option,options,or,proc,quit,read,return,save,stop,subset,then,% 
		to,try,union,use,uses,while,xor},% 
	sensitive=true,% 
	morecomment=[l]\#,% 
	morestring=[b]",% 
	morestring=[d]"% 
}[keywords,comments,strings]% 

\DeclareRobustCommand{\maketitle}{%
	\thispagestyle{fancyfirst}
	\begingroup
	\begin{center}
	\Large \textbf{Automated Counting of Restricted Motzkin Paths}\\[\medskipamount]
	\large AJ Bu \\[\medskipamount]
		\vspace{5 pt}
 \end{center}
	\endgroup
}


\fancypagestyle{fancyfirst}{%
  \pagestyle{fancy}%
  \chead{\small Automated Counting of Restricted Motzkin Paths}%
}

\begin{document}
\maketitle
 \begin{abstract} 
 Motzkin paths are simple yet important combinatorial objects. In this paper, we consider families of Motzkin paths with restrictions on peak heights, valley heights, upward-run lengths, downward-run lengths, and flat-run lengths. This paper presents two fully automated methods for enumerating the paths of such families. The first method uses numeric dynamic programming. While this method often times works, it can be slow and may not work for larger problems. The second method uses symbolic dynamic programming to solve such problems. These methods are implemented in the maple packages accompanying this article.

  \end{abstract}

\section{Introduction}%
\label{sec:introduction}

Doron Zeilberger introduced methods of counting restricted Dyck paths using numeric dynamic programming and symbolic dynamic programming in his paper "Automatic Counting of Restricted Dyck Paths via (Numeric and Symbolic) Dynamic Programming" \cite{z}. Here, I generalize his findings to the Motzkin paths. This paper is accompanied by two maple packages, which are Motzkin analogues to Zeilberger's maple packages in \cite{z}.

\begin{Def}
A \textit{Motzkin path of length $n$} is a walk in the $xy-$plane from the origin $(0,0)$ to $(n,0)$ with atomic steps $U:=(1,1)$, $D:=(1,-1)$, and $F:=(1,0)$ that never goes below the $x$-axis.
\end{Def}

For example, the following paths are Motzkin paths of length 6:
\[UUUDDD, UDUFDF, UUFFDD, FUDFUD, FFFFFF.  \]

To count the number of Motzkin paths with a given set of restrictions, let $\mathcal{P}$ denote the set of such paths and consider the weight enumerator  
  \[f(x):=\sum_{P \in \mathcal{P}} x^{Length(P)}.\]
  Note that this equals the ordinary generating function
  \[\sum_{n=0}^\infty a(n)x^n,\]
of the sequence $\{a(n)\}_{n=0}^\infty$, counting the Motzkin paths of length $n$ with the desired restrictions.

This paper presents two methods for finding the polynomial $F(x,P)$ that is zero when $P:=\sum_{n=0}^\infty a(n)x^n$.
For example, let $\mathcal{P}$ denote the set of all Motzkin paths. Note that $P\in \mathcal{P}$ either is the empty path, begins with the step $F$, or begins with the step $U$. If $P$ begins with the step $F$, then we can write 
\[P=F P_0, \]
and it is obvious that $P_0$ must also be a Motzkin path. If $P$ begins with the step $U$, then let $D_0$ denote the first time $P$ returns to the $x-$axis and write
\[P=U P_1 D_0 P_2.\]
It is easy to see that $P_1$ must be a Motzkin path shifted to height $1$ and $P_2$ is also a Motzkin path. Note that, for the paths in $\mathcal{P}$, these decompositions are unambiguous.
Moreover, given any Motzkin paths $P_0,P_1,$ and $P_2$, it is clear that the empty path, $FP_0$, and $UP_1DP_2$ are also Motzkin paths.
$\mathcal{P}$ therefore has the \emph{grammar}
\[\mathcal{P}=\{EmptyPath\} \cup F \mathcal{P} \cup U \mathcal{P} D \mathcal{P}.\]
Thus, setting $P$ equal to the weight enumerator of $\mathcal{P}$, we get the recurrence
\[P=1+xP+x^2P^2.\]

There are a fair number of papers that discuss the enumeration of certain families of Motzkin paths -- \cite{ELY}, \cite{BPPR}, and \cite{ABBG} to name a few. Recall that Dyck paths are also a family of restricted Motzkin paths, as they are Motzkin paths with no flat steps.  In "Automatic Counting of Restricted Dyck Paths via (Numeric and Symbolic) Dynamic Programming" \cite{z}, Zeilberger considers Dyck paths with restrictions on peak heights, valley heights, upward-runs, and down-ward runs. In this paper, we will look at similar restrictions.  Due to the allowance of flat-steps in Motzkin paths, however, we reevaluate what peaks and valleys are. We also introduce restrictions on flat-runs.

Given a sequence of steps $L$, define $L^n$ to be the repetition of $L$ $n$ times. (For example, $F^2 = FF$ and $(UD)^3 = UDUDUD$.) The restrictions we will consider are defined as follows:

\begin{Def} 
A \emph{peak} on a  Motzkin path is the sequence of steps $UF^kD$ for $k \geq 0$. The \emph{height} of this peak is given by the $y-$coordinate of the Motzkin path after the step $U$.
\end{Def}

\begin{Def} 
 A \emph{valley} on a Motzkin path is the sequence of steps  $DF^kU$ for $k \geq 0$. Its height is given by the $y-$coordinate after the step $D$.
\end{Def}

\begin{Def}
 A Motzkin path has an \emph{upward-run of length $n$} if it contains a run $U^n$ that is not followed by nor directly follows an up-step. 
\end{Def}

\begin{Def}
 A Motzkin path has a \emph{downward-run of length $n$} if it contains a run $D^n$ that is not directly followed by nor directly follows a down-step. 
\end{Def}
\begin{Def}
  A Motzkin path has a \emph{flat-run of length $n$} if it contains a run $F^n$ that is not directly followed by nor directly follows a flat-step. 
\end{Def}

\section{The Maple Packages}
\label{sec:package}
This article is accompanied by the following maple packages:
\begin{itemize}
    \item \texttt{Motzkin.txt}: Uses numeric dynamic programming to generate sufficiently many terms of the sequence of Motzkin paths with the desired restrictions, and then guesses the recurrence to get the desired equation.
    \item \texttt{MotzkinClever.txt}: Generates a finite system of algebraic equations by using symbolic dynamic programming and then solves the system to get the equation satisfied by the generating function directly.
\end{itemize}

These packages, example input and output files, and this article can all be found at
\\\texttt{\href{https://sites.math.rutgers.edu/~ab1854/Papers/AutocountMotzkin/AutocountMotzkin.html}{https://sites.math.rutgers.edu/\textasciitilde ab1854/Papers/AutocountMotzkin/AutocountMotzkin.html}}.

\section{Numeric Dynamic Programming  (\texttt{Motzkin.txt})}
Let us start by looking at the most basic case - finding the number of all Motzkin paths of length $N$. By definition, every Motzkin path must end with either a down-step or a flat-step. If a Motzkin path ends with a downwards-run on length $r$, then the preceding run is either an ascending-run or a flat-run that ends at height $r$. We introduce the following notation.
\begin{align*}
    u(m,n) = & \text{ the number of walks from $(0,0)$ to $(m,n)$ that never goes below the $x-$axis} \\ & \text{ and ends with an up-step.}\\
    d(m,n) = & \text{ the number of walks from $(0,0)$ to $(m,n)$ that never goes below the $x-$axis} \\ & \text{ and ends with a down-step.}\\
    f(m,n) = & \text{ the number of walks from $(0,0)$ to $(m,n)$ that never goes below the $x-$axis} \\ & \text{ and ends with a flat-step.}
\end{align*}

These give us the following equalities:

\begin{align*}
    d(m,n) & = \sum_{r=1}^m u(m-r,n+r) + f(m-r,n+r),\\
f(m,n) & =\sum_{r=1}^m u(m-r,n) + d(m-r,n), \text{ and}\\
u(m,n) & = \sum_{r=1}^m f(m-r,n-r) + d(m-r,n-r),
\end{align*}
with the initial conditions $f(0,0)=0=u(0,0)$ and $d(0,0)=1$, and the boundary conditions $d(m,k)=u(m,k)=f(m,k)=0$ for $k>m$.

\texttt{Motzkin.txt} implements these equations through the procedures \texttt{u(m,n)}, \texttt{d(m,n)}, and \texttt{f(m,n)}. Thus, to get the first $N+1$ terms of the sequence $\{a(n)\}_{n=0}^\infty$ where $a(n)$ is defined to be the number of Motzkin paths of length $n$, run
\[\texttt{seq(d(m,0)+f(m,0),m=0..N)}\]
For example,
\[\texttt{seq(d(m,0)+f(m,0),m=0..10)}\]
outputs
\[\texttt{1, 1, 2, 4, 9, 21, 51, 127, 323, 835, 2188.}\]

\subsection{Restricted Motzkin Paths}
Let $A,B,C,D$ and $E$ be arbitrary sets of positive integers -- either finite sets or infinite sets defined by the union of arithmetic progressions. We consider restricted Motzkin paths that avoid
\begin{itemize}\setlength{\itemindent}{.5in}
\item peak heights in $A$,
\item valley heights in $B$,
\item upward-runs with lengths in $C$, 
\item downward-runs with lengths in $D$, and
\item flat-runs with lengths in $E$.
\end{itemize}

In coming up with an analogue to $u(m,n)$, $d(m,n)$, and $f(m,n)$, we notice that flat-runs complicate how we count paths with restrictions on peak heights and valley heights. For example, the path 
\[UUFUDDD\]
avoids peaks with height $2$ even though it contains an upward-run that ends at height $2$. To address this, we need to define subcases for $u(m,n), d(m,n)$ and $f(m,n)$ as follows:

\begin{align*}
u_d(m,n) & = \begin{cases} 
      0 & \text{if } n \in A \\
      u(m,n) & \text{otherwise}
   \end{cases}\hspace{10pt} \text{,}\\
d_u(m,n) & =\begin{cases} 
      0 & \text{if } n \in B \\
      d(m,n) & \text{otherwise}
   \end{cases}\hspace{10pt},\\
    f_u(m,n)& = \sum_{\underset{r \not \in E}{1 \leq r \leq m}} u(m-r,n) + d_u(m-r,n)), \text{ and}\\
f_d(m,n) & =\sum_{\underset{r \not \in E}{1 \leq r \leq m}} u_d(m-r,n) + d(m-r,n).
\end{align*}

When counting restricted paths from $(0,0)$ to $(m,n)$ ending in a downward run of length $r$, the preceding run is either an upward-run or a flat-run. If it is preceded by an upward-run that ends at a height in $A$, then the path violates the restriction on peak heights. Thus, we only want to consider the paths counted by $u_d(m-r,n+r)$. Otherwise, it is preceded by a flat-run. If this flat-run is preceded by an upward-run ending at a height in $A$, then the path again has a forbidden peak height. Thus, we are interested in exactly the paths counted by $f_d(m-r,n+r)$.  Similarly, when counting restricted paths from $(0,0)$ to $(m,n)$ ending in an upward-run of length $r$, we only consider the paths counted by $d_u(m-r,n-r)$ and $f_u(m-r,n-r)$ to avoid forbidden valley heights. Note that our definitions of $f_u(m,n)$ and $f_d(m,n)$ ensure that the sub-path being counted does not end in a flat-run of length in $E$. We can use similar restrictions to ensure that our paths do not contain any forbidden run lengths.

We set 
\begin{align*}
    d(m,n) & = \sum_{\underset{r \not \in D}{1 \leq r \leq m}} u_d(m-r,n+r) + f_d(m-r,n+r),\\
f(m,n) & =\sum_{\underset{r \not \in E}{1 \leq r \leq m}} d(m-r,n) + u(m-r,n),\text{ and}\\
u(m,n) & =\sum_{\underset{r \not \in C}{1 \leq r \leq m}} d_u(m-r,n-r) + f_u(m-r,n-r). 
\end{align*}

These functions are implemented in \texttt{Motzkin.txt} and are used to get
\[\texttt{SeqABCDE(A,B,C,D,E,N)} \text{ and } \texttt{SeqABCDEr(A,B,C,D,E,r,N)},\]
which generate the terms $a(n)$ -- the number of Motzkin paths of length $n$ with the desired restrictions -- for $0\leq n \leq N$. \texttt{SeqABCDE(A,B,C,D,E,N)} is used when $A,B,C,D$ and $E$ are finite sets of non-negative integers, and \texttt{SeqABCDEr(A,B,C,D,E,r,N)} is used when the sets are defined by linear equations.

For example,
\[\texttt{SeqABCDE(\{\},\{\},\{1\},\{1\},\{1\},11)}\]
outputs
\[\texttt{[1, 0, 1, 1, 2, 1, 5, 4, 12, 13, 34, 38]},\]
and
\[\texttt{SeqABCDEr(\{2*r+1\},\{2*r+1\},\{\},\{\},\{\},r,11)}\]
outputs
\[\texttt{[1, 1, 1, 1, 2, 6, 16, 36, 73, 145, 301, 661]}.\]

The first output tells us, for example, that there are four Motzkin paths of length $7$ avoiding upward, downward, and flat runs of length $1$. We can verify that this is true by noting that such paths must either be all flat steps or a permutation of three consecutive flat-steps, two consecutive up-steps, and two consecutive down-steps. Since the up-steps must occur before the down-steps by the definition of Motzkin paths, the set of desired paths is
\[\{FFFFFFF,FFFUUDD,UUDDFFF,UUFFFDD\}.\]
The second output states that there are six Motzkin paths of length $5$ avoiding peaks and valleys with odd heights. We can easily check that the set of such paths is 
\[\{FFFFF, FUUDD, UFUDD, UUDDF, UUDFD, UUFDD\} .\] 

Note that 
\[\texttt{SeqABCDEr(\{\},\{\},\{\},\{\},\{r+1\},r,30)}.\]
outputs the number of Motzkin paths of length $n$ avoiding flat-steps for $n=0,\dots,30$.
This outputs $0$ when $n$ is odd, and the terms for even $n$ give us the list
\[\texttt{1,1,2,5,14,42,132,429,1430,4862,16796,58786,208012,742900,2674440,9694845.}\]
Inputting these terms into OEIS, we can easily verify that this is in fact the sequence of the number of Dyck paths of semi-length $n$.

\subsection{Finding the Equation Satisfied by the Generating Function}
The desired $F(x,P)$ is a polynomial, so there exist polynomials $q_0(x),\dots,q_d(x)$ such that
\[F(x,P)=q_0(x)+q_1(x)P+...+q_d(x)P^d.\] $F(x,P)$ is zero when $P:=f(x)$, the generating function of the desired sequence, thus $f(x)$ is algebraic.  Therefore, $f(x)$ satisfies a linear differential equation with polynomial coefficients, and so there is a linear recurrence equation with polynomial coefficients for the terms $a(n)$ in our sequence. (For more details see \cite{KP}, particularly Sections 6.2 and 7.2.) To get the desired polynomial, we borrow directly from Zeilberger's method of using undetermined coefficients to guess the recurrence used in \texttt{Dyck.txt} in \cite{z}.

\section{Symbolic Dynamic Programming (\texttt{MotzkinClever.txt})}

\texttt{MotzkinClever.txt} uses symbolic dynamic programming to find $F(x,P)$. More specifically, the recurrence for the set of restricted Motzkin paths is expressed as a polynomial by assigning different variables to different sets of restrictions. In addition to our original set of restricted Motzkin paths, we look at the "children" of this set. These are sets of Motzkin paths with other restrictions such that any element of our original set can be written in some form concatenating certain steps with such paths. This process is described more concretely below. We then continue to look at the children of the new sets until no new children can be produced. We will see that this must happen eventually, yielding a finite system of polynomial equations that contains the same number of equations as variables.  We use this system of equations to find the equation satisfied by the generating function.

\subsection{Avoiding Peak and Valley Heights in Finite Sets: \texttt{fAB(A,B,x,P)}}
Let $A$ and $B$ be two arbitrary finite sets of non-negative integers. We consider the ordinary generating function $f_{A,B}$ of the sequence of Motzkin paths avoiding
\begin{itemize}\setlength{\itemindent}{.5in}
    \item peak-heights in A, and
    \item valley-heights in B.
\end{itemize}
First, note that the sequence of walks with only flat-steps with weight $x^{Length(P)}$ has the generating function $\sum_{n=0}^\infty x^n = \frac{1}{1-x}.$ For convention, we say that a path has a peak at height $0$ if and only if it is a flat run. 
Now, let $\mathcal{P}$ denote the set of Motzkin paths avoiding peak heights in $A$ and valley heights in $B$, and let $\mathcal{F}$ denote the set of flat runs. Consider the following three cases:

\begin{enumerate}
    \item[Case 1:] If $0 \in A$ then let $A_1:= A \backslash \{0\}$.
    
    Let $\mathcal{P}_1$ be the set of Motzkin paths avoiding peak heights in $A_1$ and valley heights in $B$. Then it is clear $\mathcal{P}$ is the union of the disjoint sets $\mathcal{F}$ and $\mathcal{P}$, giving the following grammar
    \[\mathcal{P} \cup \mathcal{F}= \mathcal{P}_1.\]
    This gives us the following equation
    \[f_{A,B}(x)= f_{A_1,B}(x)-\frac{1}{1-x}.\]

    \item[Case 2:] If $0\not\in A$ and $0\in B$ then let $A_1:= \{a-1 | a \in A\}$ and $B_1:=\{b-1| b\in B\backslash \{0\}\}$. 
    
    Let $\mathcal{P}_1$ denote the set of Motzkin paths avoiding peak heights in $A_1$ and valley heights in $B_1$. Then any non-flat path $P$ in $\mathcal{P}$ starts with either an up-step or a flat-run followed by an up-step, and ends with either a down-step or a down-step followed by a flat-run. Note that, since $P$ avoids valleys with height $0$, it can only return to the $x-$axis once. We can therefore write
    \[P= F^{k_1} U P_1 D F^{k_2}, \]
    where $k_1$ and $k_2$ are non-negative integers, and $P_1$ is some path in $\mathcal{P}_1$. Thus, we get the grammar
    \[\mathcal{P} = \mathcal{F} \cup \mathcal{F} U \mathcal{P}_1 D \mathcal{F},\]
    which gives the following equation
  \[f_{A,B}(x)= \frac{1}{1-x}+\frac{x^2}{(1-x)^2} f_{A_1,B_1}(x).\]
    
    \item[Case 3:] If $0 \not \in A$ and $0 \not \in B$ then let  $A_1:= \{a-1 | a \in A\}$ and $B_1:=\{b-1| b\in B\}$. 
    
    Let $\mathcal{P}_1$ denote the set of Motzkin paths avoiding peak heights in $A_1$ and valley heights in $B_1$. Any non-flat path $P$ in $\mathcal{P}$ starts with either an up-step or a flat-run followed by an up-step. Then, letting $D$ denote the first time $P$ returns to the $x-$axis, we can write
    \[P= F^k U P_1 D P',\]
 where $k$ is a non-negative integer, $P_1$ is some path in $\mathcal{P}_1$, and $P'$ some path in $\mathcal{P}$. We then have the grammar
 \[\mathcal{P}= \mathcal{F} \cup \mathcal{F} U \mathcal{P}_1 D \mathcal{P}.\]
 Hence,
    \[f_{A,B}(x)=\frac{1}{1-x} + \frac{x^2}{1-x}f_{A,B}(x) f_{A_1,B_1}(x).\]
\end{enumerate}
  Thus, $\mathcal{P}$ has the child $\mathcal{P}_1.$ We then apply this procedure to $\mathcal{P}_1$ and so on. Note that we will eventually remove all the elements of $A$ and $B$ and will therefore have finitely many "descendants" of our original set. Moreover, since we have an equation to find the children of each variable, we have as many equations as variables. Each equation has only two variables, except the last equation which has one, and the variables are raised to degree at most 1. Thus, we can eliminate every variable except the one representing our original $f_{A,B}$ from the first generated equation. This gives us the polynomial satisfied by the generating function of the Motzkin paths with the desired restrictions.
  
  This procedure is implemented in \texttt{MotzkinClever.txt} by the procedure \texttt{fAB(A,B,x,P)}. For example, say we want the equation satisfied by the generating function of the sequence $\{a(n)\}_{n=0}^\infty$, where $a(n)$ is the number of Motzkin paths avoiding peak heights in $\{1,4\}$ and valley heights in $\{1,3\}$. Running
  \[ \texttt{fAB(\{1,4\},\{1,3\},x,P)}\]
 outputs the polynomial
    \begin{center} $x^8 - 2x^7 + 5x^6 - 12x^5 + 29x^4 - 38x^3 + 25x^2 - 8x + 1 + (x^6 - 16x^3 + 24x^2 - 12x + 2)(-1 + x)^3P 
      + (x^6 + 2x^5 - x^4 - 8x^3 + 12x^2 - 6x + 1)(-1 + x)^4P^2.$\end{center}
    Setting this polynomial equal to zero gives us the desired equation.
    
\subsection{Avoiding Peak and Valley Heights in Infinite Sets: \texttt{fABr(A,B,r,x,P)}}
Let $A$ and $B$ be two sets of arithmetic progressions $ar+b$ for non-negative integers $a$ and $b$. Slight modifications to the procedure \texttt{fAB(A,B,x,P)} give us the procedure \texttt{fABr(A,B,r,x,P)}, which outputs the polynomial $F(x,P)$ such that $F(x,P)=0$ is satisfied by the generating function for the sequence of Motzkin paths avoiding peak heights in $A$ and valley heights in $B$.

For example,
\[\texttt{fABr(\{2*r+1\},\{2*r+1\},r,x,P)}\]
outputs
\begin{center} $(-1 + x)^2 + (-1 + x)^3P + x^4P^2$. \end{center}

Thus, the generating function of the sequence enumerating the Motzkin paths avoiding odd peak and valley heights satisfies the equation
\[(-1 + x)^2 + (-1 + x)^3P + x^4P^2=0.\]

\subsection{Avoiding Upward-Run Lengths, Downward-Run Lengths, and Flat-Run Lengths in Finite Sets: \texttt{fCDE(C,D,E,x,P)}}
Let $C,D,$ and $E$ be finite sets of positive non-negative integers. Here, we want to find the generating function $f_{C,D,E}$ of the sequence of Motzkin paths avoiding
\begin{itemize}\setlength{\itemindent}{.5in}
    \item upward-runs with lengths in $C$,
    \item downward-runs with lengths in $D$, and
    \item flat-runs with lengths in $E$.
\end{itemize}

Let $h_{C,C_1,D,D_1,E,E_1,E_2}(x)$ weight enumerate Motzkin paths such that
\begin{itemize}\setlength{\itemindent}{.5in}
    \item the initial run is not an upward-run with length in $C_1$ nor a flat-run with length in $E_1$,
    \item the initial run is an upward-run if $0\in C_1$, and a flat-run if $0\in E_1$,
    \item the final run is not a downward run with length in $D_1$ nor a flat-run with length in $E_2$,
    \item the final run is a downward-run if $0\in D_1$, and a flat-run if $0\in E_2$,
    \item all remaining upward-run lengths are not in $C$,
    \item all remaining downward-run lengths are not in $D$, and
    \item all remaining flat-run lengths are not in $E$.
\end{itemize}

Let $\mathcal{P}$ denote the set of such paths.  For any path $P$ in $\mathcal{P}$, $P$ either leaves the $x-$axis no more than once or it can be uniquely written as 
\[P=P_1 P_2 P_3,\]
where 
\begin{enumerate}
    \item[-] $P_1$ is a path that leaves the $x-$axis no more than once and has the same restrictions as paths in $\mathcal{P}$ except it ends in a downward-run avoiding lengths in $D$,
    \item[-] $P_2$ is a path avoiding upward-runs with lengths in $C$, downward-runs with lengths in $D$, and flat-runs with lengths in $E$, and
    \item[-] $P_3$ is a path that leaves the $x-$axis no more than once and has the same restrictions defined in $\mathcal{P}$ except it begins with an upward run avoiding lengths in $C$.
\end{enumerate}

Let $H_{C,C_1,D,D_1,E,E_1,E_2}(x)$ enumerate the Motzkin paths counted by $h_{C,C_1,D,D_1,E,E_1,E_2}(x)$ that leave the $x-$axis no more than once.Then we have
\begin{align*}
       h_{C,C_1,D,D_1,E,E_1,E_2}(x)= & H_{C,C_1,D,D_1,E,E_1,E_2}(x) \\
       &\indent + H_{C,C_1,D,D \cup \{0\},E,E_1,E_2}(x)h_{C,C,D,D,E,E,E}H_{C,C \cup \{0\},D,D_1 ,E,E_1,E_2}(x).
\end{align*}
 
Note that we will never have $0 \in C_1$ and $0 \in E_1$ or $0 \in D_1$ and $0 \in E_2$, since the former statement says that the path starts with both an up-step and a flat-step, and the latter states that the path ends with both a down-step and a flat-step. To get the desired system of equations and find the children of the set $\mathcal{P}$ of paths weight-counted by  $H_{C,C_1,D,D_1,E,E_1,E_2}$, let $P$ be any path in $\mathcal{P}$ and consider the following cases.
 \begin{enumerate}
 
     \item[Case 1:] If $0 \in E_1$, then $P$ begins with a flat-step. Let $E_1'=\{e-1| e \in E_1 \backslash \{0\} \}$. Then we can write 
     \[P=F P_1,\]
where $P_1$ is a path weight-counted by $H_{C,C,D,D_1,E,E_1',E_2}(x)$. Hence, 
     \[H_{C,C_1,D,D_1,E,E_1,E_2}(x)=xH_{C,C,D,D_1,E,E_1',E_2}(x)\]
     
     \item[Case 2:] If  $0 \not \in E_1$  and $0 \in E_2$ then $P$ ends with a flat-step. Let $E_2'=\{e-1| e \in E_2 \backslash\{0\} \}$, and write
     \[ P= P_1 F,\]
     where $P_1$ is a path weight-counted by $H_{C,C,D,D_1,E,E_1,E_2'}(x)$. This gives us
     \[H_{C,C_1,D,D_1,E,E_1,E_2}(x)=xH_{C,C,D,D_1,E,E_1,E_2'}(x)\]
     
     \item[Case 3:] If $0 \not \in E_1$, $0 \not \in E_2$, $0 \in C_1$, and $0 \in D_1$, then $P$ starts with an up-step and ends with a down-step. Letting
     $C_1'=\{c-1| c \in C_1 \backslash\{0\} \}$ and
     $D_1'=\{d-1| d \in D_1 \backslash\{0\}\}$, we can write
     \[P=U P_1 D,\]
     where $P_1$ is a path weight-counted by $h_{C,C_1',D,D_1',E,E_1,E_2}(x)$. (Note that $P_1$ is able to return to the height it begins at more than once.) Thus,
     \[H_{C,C_1,D,D_1,E,E_1,E_2}(x)=x^2h_{C,C_1',D,D_1',E,E,E}(x)\]
     
     \item[Case 4:] If $0 \not \in E_1$, $0 \not \in E_2$, $0 \not \in C_1$, and $0 \not \in D_1$ then $P$ is either the empty path, starts with an up-step, or starts with a flat-step. We therefore get
     \[H_{C,C_1,D,D_1,E,E_1,E_2}(x)=H_{C,C_1 \cup \{0\} ,D,D_1,E,E,E_2,}(x)  + H_{C,C_1,D,D_1,E,E_1 \cup \{0\},E_2}(x)+1\]

     \item[Case 5:] If $0 \not \in E_1$, $0 \not \in E_2$, $0 \not \in C_1$, and $0 \in D_1$, then $P$ is non-empty and starts with either an up-step or a flat-step. Hence,
        \[H_{C,C_1,D,D_1,E,E_1,E_2}(x)=H_{C,C_1 \cup \{0\} ,D,D_1,E,E,E_2}(x)  + H_{C,C_1,D,D_1,E,E_1 \cup \{0\},E_2}(x)\]
     
     \item[Case 6:] If $0 \not \in E_1$, $0 \not \in E_2$, $0  \in C_1$, and $0 \not \in D_1$, then $P$ is non-empty and ends with either a down-step or a flat-step. Thus,
        \[H_{C,C_1,D,D_1,E,E_1,E_2}(x)=H_{C,C_1,D,D_1 \cup \{0\},E,E,E_2}(x) + H_{C,C_1,D,D_1,E,E_1,E_2\cup \{0\}}(x).\]
 \end{enumerate}
We again generate finitely many descendants from the original set, as we will eventually remove all of the elements in $C$, $D$, and $E$. We also have as many equations as variables. Note that these polynomials generate an ideal. Since any basis will give the same set of solutions, we can look at the reduced Gr\"{o}bner basis of the generated ideal.

\subsubsection{Gr\"{o}bner Bases - A Quick Background and Their Application}

\begin{Def}
A \emph{Gr\"{o}bner basis} of an ideal $I \subset k[x_1,...,x_n]$ (with respect to a monomial order >) is a finite subset $G=\{g_1,...,g_t\}$ of $I$ such that for that every nonzero polynomial $f$ in $I$, the leading term of $f$ is divisible by the leading term of $g_i$ for some $i$.
\end{Def}

\begin{Def}
A Gr\"{o}bner basis $G$ is \emph{reduced} if, for every element $g \in G$, no monomial in $g$ is in $\la LT(G-{g})\ra,$ the ideal generated by the leading terms of the other elements in $G$.
\end{Def}

Choosing the correct monomial ordering (namely, pure lexicographic order) will allow us to ensure that the smallest element of the reduced Gr\"{o}bner basis is in the form to most easily find the desired equation satisfied by the generating function. This is due to the following theorem. (Here, we assign each descendant found in our system of equations a variable $x_i$, and let $x_n$ be the variable representing the original family of restricted Motzkin paths. We do not consider $x$ as one of these variables.)

\begin{Thm}[The Elimination Theorem]
If $G$ is a Gr\"{o}bner basis for $I$ with respect to lex order $x_1>x_2 > \dots > x_n$, then \[G_\ell = G \cap k[x_{\ell+1},...,x_n]\]
is a Gr\"{o}bner basis of the $\ell$-th elimination ideal $I_\ell=I \cap k[x_{\ell+1},...,x_n]$.
\end{Thm}

If $f(x)$ denotes the generating function of the sequence enumerating our original family of restricted Motzkin paths, then $x_n=f(x)$ is a partial solution to our system of equations represented by $I$. By the Elimination Theorem, if $q$ denotes the smallest polynomial of the reduced Gr\"{o}bner basis, then either $q \in G_{n-1}$ or $I_{n-1}=\la 0 \ra$.  $I_{n-1}= \la 0 \ra$, however, contradicts the existence of the desired nonzero polynomial $F(x,P).$ Thus, $q \in G$ is a polynomial in terms of $x$ and $x_n$ and is zero when $x_n=f(x).$  Factoring $q$ completely, we can write 
\[q=q_1^{d_1}\dots q_k^{d_k},\]
where $d_i \geq 1$. If $k=1$, then we are done and $F(x,P):=q_1$. Otherwise, $x_n=f(x)$ also satisfies $q_i=0$ for one of the factors $q_i.$  We can then use the first $m$ terms, where $m$ is sufficiently large, of the sequence of interest to determine which factor is the desired $q_i$. We thereby get the desired polynomial $F(x,P):=q_i.$

This process is implemented in $\texttt{fCDE(C,D,E,x,P)}.$
 \[\texttt{fCDE(\{1,2,3\},\{\},\{\},x,P)}\]
outputs
\[1 + (-x^2 + x - 1)P - x^2(x - 1)P^2 + P^4x^8 + P^5x^9,\]
and
 \[\texttt{fCDE(\{\},\{1\},\{1\},x,P)}\]
 outputs
\[x^2 - x + 1 + (-x^4 + x^3 - x^2 + x - 1)P + x^2(x^4 - x^3 + x^2 - x + 1)P^2 + P^3x^6.\]

Thus, when $P=\sum_{n=0}^\infty a(n)x^n,$ where $a(n)$ is the number of Motzkin paths of length $n$ avoiding upward runs of lengths 1, 2, and 3,
\[1 + (-x^2 + x - 1)P - x^2(x - 1)P^2 + P^4x^8 + P^5x^9=0.\]
If $a(n)$ is the number of Motzkin paths of length $n$ avoiding downward-runs and flat-runs of length 1, then
\[x^2 - x + 1 + (-x^4 + x^3 - x^2 + x - 1)P + x^2(x^4 - x^3 + x^2 - x + 1)P^2 + P^3x^6=0.\]

\subsection{Avoiding Upward-Run Lengths, Downward-Run Lengths, and Flat-Run Lengths in Infinite Sets: \texttt{fCDEr(C,D,E,r,x,P)}}
Suppose $C,D,$ and $E$ are sets of arithmetic progressions $ar+b$ for non-negative integers $a$ and $b$. Through slight modifications to \texttt{fCDE(C,D,E,x,P)}, we get the procedure \texttt{fCDEr(C,D,E,r,x,P)}. \texttt{fCDEr(C,D,E,r,x,P)} finds the desired polynomial $F(x,P)$ that is zero when $P$ is the generating function for the sequence of Motzkin paths avoiding upward-run lengths in $C$, downward-run lengths in $D$, and flat-run lengths in $E$.
Running
\[\texttt{fCDEr(\{2*r+1\},\{2*r+1\},\{2*r+1\},r,x,P)}\]
tells us that when $P$ is the generating function of the sequence enumerating Motzkin paths avoiding upward, downward, and flat runs of odd length, we have
\[1 + (x - 1)(x + 1)P + P^2x^4=0\]

To get the equation satisfed by the generating function of the sequence enumerating Motzkin pats avoiding upward runs of odd lengths and flat-runs of positive even length, input
\[\texttt{fCDEr(\{2*r+1\},\{\},\{2*r+2\},r,x,P)}.\]
This tells us that our desired equation is
\[x^2 - x - 1 - (x - 1)(x + 1)P + x^4(x^2 - x - 1)P^3=0.\]

\section{Conclusion}
Using similar approaches, we can create ways to automate counting of other objects. The approach of using numeric dynamic programming can efficiently generate many terms of the desired sequence. Guessing the algebraic equation, however, will not always work well. Thus, for larger problems, we need to use symbolic dynamic programming instead. Here, we identify recursive relations for the sets of relevant objects. Then, we use this system of equations to find an equality solved by the weight-enumerator of the set of combinatorial objects of interest.

%\makefoot
\end{document}